\documentstyle[12pt,amscd,amssymb]{amsart}

\newtheorem{th}{Theorem}[section] 
\newtheorem{lem}[th]{Lemma} 
\newtheorem{cor}[th]{Corollary} 

\theoremstyle{definition} 
\newtheorem{dfn}[th]{Definition} 
 
 \newtheorem{question}[th]{Question}

\newtheorem{ass}[th]{Assumption}

\theoremstyle{remark} 
\newtheorem{rem}[th]{Remark}

\newcommand{\rar}{\rightarrow} 
 
\newcommand{\bfp}{{\Bbb{P}}} 
\newcommand{\bfq}{{\Bbb{Q}}} 
 
\newcommand{\bfr}{{\Bbb{R}}} 
 
\newcommand{\co}{{\cal{O}}}
\newcommand{\cm}{{\cal{M}}}
  
 \newcommand{\setmin}{\,{^{_\setminus}}\,}
\newcommand{\QED}{\qed\\}

\newcommand{\spec}{{\operatorname{Spec\ }}}

\newcommand{\dimension}{{\operatorname{dim\ }}}

\newcommand{\proj}{{\operatorname{\bf Proj\ }}}

\begin{document}

\noindent%

\title[Boundedness]{Minimal models and boundedness of stable varieties\\
Preliminary version}
\author{Kalle Karu}
\address{Department of Mathematics\\ Boston University\\ 111 Cummington\\
Boston, MA 02215, USA} 
\email{kllkr@@math.bu.edu}
\date{\today}

\maketitle

\section{Introduction}

It is known that there exists a quasi-projective coarse moduli space for
nonsingular n-dimensional varieties $X$ with ample canonical divisor $K_X$
and given Hilbert function $H(l) = \chi(X, \co(l K_X))$ (e.g. see \cite{vie}). The
starting point in the construction of this moduli space is Matsusaka's big
theorem \cite{mat1}  which states that the class of these varieties is
bounded: there exists an integer $\nu_0 >0$, independent of $X$, such that
$\nu_0 K_X$ is 
very ample and has no higher cohomology. This condition is equivalent to the
existence of a flat family $\cal{X}\rar B$ over a scheme of finite type whose
fibers contain all such varieties.

In order to compactify the moduli space of smooth n-folds, one has to add
points at the boundary, corresponding to limits of smooth varieties. The natural
way to construct these limits is to take a one-parameter  semistable family of
n-folds, and find the relative canonical model of this family in the sense of
the 
minimal model program in dimension $n+1$ $(MMP(n+1))$. All such limits have
semi-log-canonical singularities. 

We consider the class of stable (almost) smoothable n-folds $X$ with
semi-log-canonical singularities and ample canonical divisor
$K_X$. Smoothability here means that $X$ 
admits a deformation to a variety with rational Gorenstein singularities. The
moduli functor $\cm_H^{sm}$ then assigns to a scheme $S$ the set of
isomorphism classes of families of stable smoothable n-folds over $S$ with a
given Hilbert function $H$. The main result of this paper is that the minimal
model program in dimension $n+1$ implies boundedness of the functor
$\cm_H^{sm}$.

\begin{th}\label{thm-bound} Assuming MMP(n+1), there exists a family
$\tilde{X}\rar B$ over a projective scheme $B$ in $\cm_H^{sm}(B)$  whose
geometric fibers include all stable smoothable n-folds with Hilbert function
$H$. 
\end{th}

If, in addition to boundedness, $\cm_H^{sm}$ is locally closed, separated,
complete, with tame
automorphisms, and if the canonical polarization is semi-positive, then
Koll\'{a}r's theorem \cite{kol} states that there exists a projective moduli
space $M_H^{sm}$ coarsely representing the functor
$\cm_H^{sm}$. Semi-positivity of the canonical polarization has been proved
by  Koll\'{a}r (Theorem 4.12, \cite{kol}). Since $K_X$ is ample, $X$ has a
finite automorphism group by a result of Iitaka (Section 11.7,
\cite{iit}). Separatedness of $\cm_H^{sm}$ follows from the uniqueness of
log-canonical models (Lemma~\ref{lem-sep}). In Lemma~\ref{lem-def-slc} we
prove that semi-log-canonical varieties deform to semi-log-canonical
varieties; since the other conditions we require from stable smoothable
n-folds are either open or closed, the functor $\cm_H^{sm}$ is locally
closed. Completeness then follows because $M_H^{sm}$ is the
image of the projective scheme $B$. This proves the existence of coarse
moduli spaces:

\begin{cor} Assuming MMP(n+1), the moduli functor $\cm_H^{sm}$ is coarsely
represented by a projective scheme $M_H^{sm}$.
\end{cor}

The minimal model program MMP(n+1) is known only in dimensions $n+1\leq
3$. Boundedness for semi-log-canonical surfaces has been proved by Alexeev
\cite{ale1}. With this proof he finished the construction of projective
coarse moduli spaces for semi-log-canonical surfaces that was started in
\cite{ksb}. The proof of Theorem~\ref{thm-bound} simplifies Alexeev's result of
boundedness for semi-log-canonical surfaces, and generalizes it to n-folds,
subject to the minimal model program assumption.

In \cite{ale2} Alexeev considered the problem of constructing projective
coarse moduli spaces for stable pairs $(X,B)$ consisting of an n-fold $X$ and
a divisor $B$ in $X$. He showed that the log-minimal model program in
dimension $n+1$ together with  the boundedness assumption imply the existence
of coarse moduli spaces for stable pairs. At the end of the paper we
indicate how the proof of Theorem~\ref{thm-bound} can be modified to prove
boundedness of stable smoothable pairs.

\subsection{Sketch of the proof.} We imitate the construction of a stable
n-fold $X$ as a limit of a  one-parameter family of (almost) smooth
varieties, and try to find all these limits at the same time. By Matsusaka's
theorem \cite{mat2}, the class of normal varieties with rational Gorenstein
singularities and ample canonical bundle is bounded. We start with a
compactification of the Hilbert scheme parameterizing these varieties, and the
universal family over it. To this morphism we apply weak semistable reduction
\cite{ak} and then take the relative canonical model of the weakly semistable
family. The crucial step in showing that this canonical model exists is
a theorem of Siu and Kawamata on the invariance of plurigenera in families of
canonical varieties \cite{kaw}. If $f:X\rar B$ is a weakly semistable
morphism, then for a general nonsingular curve $C\subset B$ the restriction
$X_C=f^{-1}(C)$ has canonical 
singularities, so we can find its relative canonical model by the $MMP(n+1)$
assumption. Now invariance of plurigenera implies that sections of
$f_*\co_{X_C}(K_{X_C})$ can be lifted to sections of $f_*\co_X(K_X)$, and it
follows from this that 
the relative canonical ring of $f:X\rar B$ is a finitely generated
$\co_B$-algebra. 

Most results we prove about semi-log-canonical varieties and the moduli
functor $\cm_H^{sm}$ are direct generalizations of the corresponding results 
for surfaces given in \cite{ksb, ale3}, or simplifications of the results for
stable pairs \cite{ale2}.

\subsection{Acknowledgement} I wish to thank J\'{a}nos Koll\'{a}r
for suggesting the problem of boundedness, and for several ideas in the
proof of it. I am also very grateful to Dan Abramovich for advice.

\section{preliminaries}

We work over an algebraically closed field $k$ of characteristic zero.

By a variety $X$ we mean a reduced (not necessarily irreducible) separated
scheme of finite type over the field $k$. A birational morphism $f:X\rar Y$
is a morphism that maps every 
component of $X$ birationally to a component of $Y$, and has a birational
inverse. Whenever we talk about a Weil divisor $D$ in $X$ we assume that the
generic point of every component of $D$ lies in the smooth locus of $X$. A
$\bfq$-linear combination of Weil divisors $D=\sum a_i D_i$  in $X$ is
$\bfq$-Cartier if a multiple $l D$ is Cartier for some integer $l>0$. If
$f:Y\rar X$ is a 
morphism, we let $f^*(D)$ be the $\bfq$-Cartier divisor $\frac{1}{l} f^*(l
D)$; and if $f$ is birational then $f_*^{-1}(D)= \frac{1}{l} f_*^{-1}(l D)$
is the strict transform of $D$.

Suppose that $X$ contains an open subvariety $U\subset X$ with complement of
codimension at least two such that $U$ is Gorenstein, i.e. the dualizing
sheaf $\omega_U$ is invertible. We choose a canonical divisor $K_U$, where  
$\co_U(K_U) \cong \omega_U$, and let $K_X$ be its closure in $X$. The variety
$X$ is called $\bfq$-Gorenstein if $K_X$ is $\bfq$-Cartier. For a morphism
$f: Y\rar X$, we let $K_{Y/X}$ be such that $\co_X(K_Y) \cong \co_Y(K_{Y/X})
\otimes f^* 
\co_X(K_X)$. 

\begin{dfn} Let $X$ be a normal $\bfq$-Gorenstein variety, and let $D$ be a
$\bfq$-Cartier divisor in $X$. We say that the pair $(X,D)$ has canonical
(resp. log
canonical) singularities if for any birational morphism $f: Y\rar X$ from a
normal $\bfq$-Gorenstein variety $Y$ we have
\[ K_{Y} = f^*(K_X+D) - f_*^{-1}(D) +\sum_i a_i E_i \]
where all $a_i \geq 0$ (resp. $a_i \geq -1$).
\end{dfn}

To show that a variety has canonical or log-canonical singularities, it
suffices to check the numerical condition for a resolution $f: Y\rar X$ of
the pair $(X,D)$, i.e. $Y$ is nonsingular and the inverse image $f^*(D)$ is a
divisor with components crossing normally.

A flat projective morphism $f: X\rar C$ from a nonsingular variety $X$ to a
nonsingular curve $C$ is semistable if its fibers are reduced divisors of
simple normal crossing.  
If $f: X\rar C$ is any flat projective morphism to a nonsingular curve, we
say  that $X$ admits a semistable resolution if there exists a birational
morphism $g: \tilde{X}\rar X$ such that $f \circ g: \tilde{X}\rar C$ is
semistable.  

To prove Theorem~\ref{thm-bound} we need to assume the existence of relative
canonical models (cf. \cite{kmm}) in the following two cases:

\begin{ass} (MMP(n+1)). Let $f: X\rar Y$ be a morphism of varieties,
$\dimension(X)=n+1$. Assume that either  
\begin{enumerate} 
\item $f$ is birational and for some morphism $g:Y\rar C$ to nonsingular
curve $C$, the composition $g\circ f: X\rar C$ is semistable; or
\item $f$ is a flat projective morphism from a variety $X$ with canonical
singularities to a nonsingular curve $Y$ with fibers of general type. 
\end{enumerate}
Then the relative canonical ring
\[ R_{X/Y} = \bigoplus_{l\geq 0} f_*\co_X(l K_X) \]
is a finitely generated $\co_Y$-algebra. The scheme $X_{can} = \proj R_{X/Y}$
(or the morphism $f_{can}: X_{can} \rar Y$) is called a relative canonical model. 
\end{ass}

Relative canonical models are unique: if $X'\rar Y$ is another
morphism satisfying conditions 1. or 2. above, and if $X'$ is birational to
$X$ over $Y$, then $X'_{can}$ is isomorphic to $X_{can}$ over $Y$. 

\subsection{Semi-log-canonical singularities}

\begin{dfn}\label{dfn-slc} (cf. \cite{ale3} Def. 2.8, \cite{kol} Def. 4.10)
We say that a $\bfq$-Gorenstein variety $X$ has semi-log-canonical
singularities if 
\begin{itemize}
\item[(i)] $X$ satisfies Serre's condition $S_2$;
\item[(ii)] $X$ has normal crossing singularities in codimension 1; 
\item[(iii)] for any birational morphism $f: Y\rar X$ from a
normal $\bfq$-Gorenstein variety $Y$ we have
\[ K_{Y} = f^*(K_X) +\sum_i a_i E_i \]
where all $a_i \geq -1$.
\end{itemize}
\end{dfn}

\begin{rem} Condition $(ii)$ implies that $X$ is Gorenstein in codimension
one, so $K_X$ is well defined. We call the closure of the singular locus in
codimension one the double divisor of $X$. Let $\nu:X^\nu\rar X$ be the
normalization and $cond(\nu)$ the reduced effective divisor defined by 
\[ K_{X^\nu} = \nu^*(K_X) - cond(\nu).\]
Then $cond(\nu)$ is the inverse image of the double divisor (with
coefficients equal to
one), and the condition $(iii)$ is equivalent to the pair
$(X,cond(\nu))$ having log-canonical singularities. Again, it suffices to check
the numerical condition for a resolution $Y$ of the pair $(X,cond(\nu))$.

The only semi-log-canonical varieties we are going to consider
are fibers of morphisms $Y\rar C$ from a variety with canonical
singularities to a nonsingular curve. These varieties are Cohen-Macaulay, so
the condition $(i)$ is automatically satisfied. If, moreover, the family admits a
semistable resolution, then $(ii)$ is also satisfied: localizing at a
codimension two point of $Y$ (codimension one in $X=Y_0$) we get a germ of a
surface with canonical singularities, admitting a semistable resolution; the
fibers of such a surface are curves with normal crossing singularities.
\end{rem}

\begin{th}\label{lem-ksb} (cf. \cite{ksb} Thm. 5.1.) Let $\pi: X\rar C$ be a flat
projective morphism from a normal $\bfq$-Gorenstein variety to a germ of a
nonsingular curve $(C,0)$. 
\begin{itemize}
\item [(i)] Suppose that $X$ admits a semistable resolution. Then $X$ has
canonical singularities iff the special fiber $X_0$ has semi-log-canonical
singularities and the general fiber $X_t$ has canonical singularities.
\item[(ii)] The product $X\times_C C'$ has canonical singularities for all
finite base changes $(C',0')\rar(C,0)$ iff the special fiber $X_0$ has
semi-log-canonical singularities and the general fiber $X_t$ has canonical
singularities. 
\end{itemize}
\end{th}

{\bf Proof.} $(i)$ $\Rightarrow$. Let $g:\tilde{X}\rar X$ be a semistable
resolution, and assume that $X$ has canonical singularities:
\[ K_{\tilde{X}} = g^* K_X +\sum_i a_i E_i +\sum_j b_j F_j,\]
where $E_i$ are the exceptional divisors mapping to $0\in C$, $F_j$ are flat
over $C$, and $a_i,b_j\geq 0$ for all $i,j$. If $X_0'=g_*^{-1} X_0$ is the strict
transform of $X_0$, then by adjunction
\[ K_{X_0'} = (K_{\tilde{X}} + X_0' )|_{X_0'} = (K_{\tilde{X}} - \sum E_i) |_{X_0'} =
(g^* K_X +\sum_i (a_i-1) E_i+\sum_j b_j F_j)|_{X_0'}. \]

Let $\nu: Y\rar X_0'$ be the normalization map. Then $K_Y = \nu^*K_{X_0'}-
cond(\nu)$ where $cond(\nu)$ is a reduced effective normal crossing divisor,
not containing any components supported on $E_i$. This proves that the
special fiber $X_0$ is semi-log-canonical because $Y\rar X_0$ is a
resolution. Since  $X$ is canonical, the general fiber $X_t$ is also
canonical. 

$\Leftarrow$. Let $g: \tilde{X}\rar X$ be the relative canonical model over
$X$ obtained by applying the MMP(n+1) assumption to a semistable resolution
of $X$. From the assumption that the general fiber of $\pi$ is canonical
it follows that $g$ is an isomorphism away from the special fibers; in
particular, there are no exceptional divisors $F_j$ flat over $C$ (with
notation as above). 

Suppose that some $a_i\geq 0$, say $a_0\geq 0$ and $a_0$ is maximal among the
$a_i$. Let $Y\subset E_0$ be a curve mapping to a point in $X$, not lying in
any other component of $\tilde{X}_0$. Then
\[ K_{\tilde{X}}.Y = \sum a_i E_i.Y = \sum_{i\neq 0} a_i E_i.Y -
\sum_{i\neq 0} a_0 E_i.Y \leq 0 \]
and this contradicts the g-ampleness of $K_{\tilde{X}}$. So, all $a_i<0$,
and since $X_0$ is semi-log-canonical, the adjunction formula above shows
that there are no exceptional divisors
$E_i$. Thus, $X$ has canonical singularities.

$(ii)$ $\Rightarrow$. We know that after a finite base change
$(C',0')\rar(C,0)$ the fiber product $X'=X\times_C C'$ admits a semistable
resolution. By part $(i)$ the special fiber $X'_{0'}\cong X_0$ is
semi-log-canonical  and the general fiber $X'_{t'}\cong X_t$ is canonical.

$\Leftarrow$. It suffices to prove that $X$ has canonical singularities. For
some finite base change $(C',0')\rar(C,0)$, the fiber product $X' = X\times_C
C'$ admits a semistable resolution and has canonical singularities by
(i). Let $h:\tilde{X}\rar X$ be a resolution of $(X,X_0)$, and let 
$\tilde{X}' = \tilde{X}\times_X X'$:
\[     
\begin{CD}
\tilde{X}' @>{\tilde{f}}>> \tilde{X}\\
@V{g}VV @VV{h}V\\
X' @>{f}>> X
\end{CD}
\]
First assume that $\tilde{X}'$ is normal. Tracing the diagram in two ways we get 
\begin{gather*}
K_{\tilde{X}'} = g^* K_{X'}+\sum a_i E_i = g^* f^* K_X + m\tilde{X}'_0 + \sum
a_i E_i,\\ 
K_{\tilde{X}'} =\tilde{f}^* K_{\tilde{X}} + m\tilde{X}'_0 = \tilde{f}^* h^* K_X
+ \tilde{f}^* \sum b_i F_i + m\tilde{X}'_0.
\end{gather*}
Since $X'$ has canonical singularities, all
$a_i\geq 0$, hence all $b_i\geq 0$ and $X$ has canonical singularities.  

If $\tilde{X}'$ is not normal, let $\nu:\tilde{X}^\nu \rar \tilde{X}'$ be the
normalization, and $cond(\nu)$ the conductor:
\[ K_{\tilde{X}^\nu} = \nu^* K_{\tilde{X}'} - cond(\nu). \]
Computing $ K_{\tilde{X}^\nu}$ in two ways as above, we get 
\[ \nu^* \tilde{f}^* \sum_i b_i F_i = \sum a_i E_i + cond(\nu),\]
where $a_i\geq 0$. Since $cond(\nu)$ is effective, we see that $b_i\geq 0$
and $X$ has canonical singularities.  \QED

The proof of the previous theorem can be easily modified to show that
semi-log-canonical singularities deform to semi-log-canonical singularities. 

\begin{lem}\label{lem-def-slc} (cf. \cite{ksb} Cor. 5.5.) Let $\pi: X\rar C$
be a flat projective morphism from a $\bfq$-Gorenstein variety
to a germ of a 
nonsingular curve $(C,0)$. If the special fiber $X_0$ has semi-log-canonical
singularities, then the general fiber $X_t$ also has semi-log-canonical
singularities.
\end{lem}

{\bf Proof.} First assume that $X$ is normal and admits a semistable
resolution. Let $g:\tilde{X}\rar X$ be the relative canonical model over $X$
of this resolution. As in the proof of part (i) of the previous theorem, write
\[ K_{\tilde{X}} = g^* K_X +\sum a_i E_i + \sum b_j F_j \]
where $E_i$ are the exceptional divisors mapping to the special fiber $X_0$,
and $F_j$ are the exceptional divisors flat over $C$.   

When we restrict to the general fiber $X_t$, it follows from the
$g$-ampleness of 
$K_{\tilde{X}}$ that all $b_j <0$ (cf. Lemma 2.19, \cite{kol1}). As before,
we get from the adjunction formula  that there are no exceptional divisors
$E_i$ mapping to $X_0$. Since the special fiber $X_0$ is semi-log-canonical,
the coefficients $b_j\geq -1$. It follows that the general fiber $X_t$ is
semi-log-canonical because $\tilde{X}_t\rar X_t$ is a resolution.

If $X$ is not normal, we let $\nu: X^\nu\rar X$ be the normalization. We may
assume that the components of $X^\nu$ admit a semistable resolution (with the
same base change) $g:\tilde{X}\rar X^\nu\rar X$.
The conductor $cond(\nu)$ maps flat onto $C$, and as in the normal case,
the discrepancies $b_j$ of the general fiber are the same as the
discrepancies of the special fiber. \QED

\subsection{The moduli functor $\cm_H^{sm}$}

\begin{dfn}\label{dfn-stable} (cf. Defs. 5.2, 5.4, \cite{kol})
\begin{enumerate}
\item
\begin{enumerate}
\item[(a)] A stable n-fold $Y_0$ is a connected projective n-dimensional
variety (not necessarily irreducible) over an algebraically closed field of
characteristic zero such that $Y_0$ has  semi-log-canonical singularities and
the canonical divisor $K_{Y_0}$ is ample.
\item[(b)] The n-fold $Y_0$ is smoothable if there exists a flat projective
$\bfq$-Gorenstein one-parameter family $\pi: Y\rar C$  of stable n-folds such
that the special 
fiber of $\pi$ is $Y_0$ and the general fiber is a normal n-fold  with at
most rational Gorenstein singularities.
\item[(c)] Given a polynomial $H$, we say that  $Y_0$ has Hilbert function
$H$ if $\chi(Y_0,\co(l K_{Y_0}))=H(l)$ for all $l\geq 0$ .
\end{enumerate}

\item A stable smoothable n-fold with given Hilbert function $H$ over a
scheme $S$  is a flat projective morphism $\pi:Y\rar S$ such that 
\begin{enumerate}
\item[(a)] every geometric fiber $Y_s=\pi^{-1}(s)$ is a stable smoothable
n-fold with Hilbert function $H$; and
\item[(b)] we have a natural isomorphism
\[ \co_{Y_s}(l K_{Y_s}) \cong \co_Y(l K_Y) |_{Y_s} \]
for every $l>0$ and every point $s\in S$.
\end{enumerate}
\end{enumerate}
\end{dfn}

\begin{dfn} We define the moduli functor $\cm_H^{sm}:
(\mbox{Schemes})\rar(\mbox{Sets})$ by 
$\cm_H^{sm}(S)$ = \{isomorphism classes of stable smoothable n-folds
with given Hilbert function $H$ over $S$\}.
\end{dfn}

Given a smoothing $\pi: Y\rar C$ of a stable n-fold $Y_0$ as in
Definition~\ref{dfn-stable}, the total space $Y$ has canonical
singularities by Theorem~\ref{lem-ksb}~(ii), and $K_Y$ is $\pi$-ample. So,
$\pi: Y\rar C$  is a relative canonical model. Conversely, if $\pi: Y\rar C$
is a relative canonical model which admits a semistable resolution, the
special fiber $Y_0$ is a stable n-fold by Theorem~\ref{lem-ksb}~(i); and if the
generic fiber of $\pi$ is Gorenstein, then $\pi: Y\rar C$ is in fact a
smoothing of $Y_0$. This shows that the moduli functor
$\cm_H^{sm}$ can be defined purely in terms of special fibers of canonical
models.

Note that the conditions $1.(b)$ and $2.(b)$ in Definition~\ref{dfn-stable}
are both closed for flat families $Y\rar S$: the first by deformation theory,
and the second by Kollar's result on push-forward and base
change. Lemma~\ref{lem-def-slc} proves that having semi-log-canonical
singularities is an open condition. The condition that $K_Y$ restricts to an
ample divisor on fibers $Y_s$ is also open (see Remark~\ref{rem-ample-open}
below). It follows that the functor  $\cm_H^{sm}$ is locally closed.

\begin{lem} \label{lem-sep} The functor $\cm_H^{sm}$ is separated; i.e.
any family $X\rar C\setmin\{0\}$ of stable smoothable n-folds over a
punctured curve can be extended to a family over $C$ in no more than one way.
\end{lem}

{\bf Proof.} (Theorem 3.28, \cite{ale3}). Suppose that we have two extensions
$f_1: X_1\rar C$ and $f_2: X_2\rar C$. Let $\tilde{X}$ be a common resolution
of $X_1$ and $X_2$, $g_1: \tilde{X}\rar X_1$, $g_2:\tilde{X}\rar X_2$. Write 
\[ K_{\tilde{X}} = g_1^* K_{X_1} +\sum a_i E_i + \sum b_j F_j \]
where $E_i$ are the exceptional divisors mapping to $0\in C$,
and $F_j$ are flat over $C$. Taking a finite base change $C'\rar C$ if
necessary, we may assume that $\tilde{X}$ admits a semistable resolution. 

First assume that $X_1$ is normal. As in the proof of Theorem~\ref{lem-ksb},
we get that $a_i\geq 0$ and $b_j\geq -1$. Since $K_{X_1}$ is
$f_1$-ample, $X_1$ is the log-canonical model of the divisor $K_{\tilde{X}}
+\sum F_j$ and is given as 
\[ X_1 \cong \proj \bigoplus_{l\geq 0} (f_1\circ g_1)_*
\co_{\tilde{X}}(l(K_{\tilde{X}} +\sum F_j)). \]
The same is true for $X_2$, and since $X_1$ and $X_2$ are isomorphic away
from the special fibers, the divisors $F_j$ are the same in both cases. Hence
$X_1$ and $X_2$ are isomorphic.

If $X_1$ is not normal, we can recover its normalization
$X^{\nu}_1$ as the log-canonical model of the divisor $K_{\tilde{X}} + \sum
F_j$, where now  $\sum F_j$ also includes the components of the inverse image
of $cond(\nu)$ flat over $C$. This proves that $X_1$ and $X_2$ are
isomorphic in codimension one. Both $X_1$ and $X_2$ satisfy Serre's condition
$S_2$, so they are in fact isomorphic. \QED

\subsection{Weak semistable reduction}

Recall from \cite{te} that $U_X\subset X$ is
a toroidal embedding if at every closed point $x\in X$ it is formally
isomorphic to a torus embedding $T_x\subset Y_x$; such $T_x\subset Y_x$ is
called a local model at $x$. A surjective morphism $f:(U_X\subset X)\rar
(U_B\subset B)$ is toroidal if it is equivariant in local models at every
point. We only consider toroidal morphisms without horizontal divisors; 
i.e. no component of $X\setmin U_X$ dominates a component of $B$.

\begin{dfn} A toroidal morphism $f:(U_X\subset X)\rar (U_B\subset B)$ without
horizontal divisors is weekly semistable if 
\begin{enumerate}
\item $f$ is equidimensional;
\item $f$ has reduced fibers; and
\item $B$ is nonsingular.
\end{enumerate}
The morphism $f$ is semistable if also $X$ is nonsingular.
\end{dfn}

The following weak semistable reduction theorem was proved in \cite{ak};

\begin{th}\label{thm-wssr} Let $X\rar B$ be a surjective morphism of
projective varieties 
with geometrically integral generic fiber. There exist a generically finite
proper surjective morphism $B'\rar B$ and a proper birational morphism
$X'\rar X\times_B B'$ such that the induced morphism $f': X'\rar B'$ is
weakly semistable. \qed
\end{th}

Using the toroidal structure of a weakly semistable morphism it is not
difficult to prove the following properties:

\begin{lem}\label{lem-wss-prop} Let $f:X\rar B$ be weakly semistable. Then
\begin{enumerate}
\item $X$ has rational Gorenstein (hence canonical) singularities.
\item If $g:C\rar B$ is a morphism from a nonsingular curve such that
$g(C)\not\subset B\setmin U_B$ then $C$ and $X_C=X\times_B C$ can be given
toroidal structures such that $X_C\rar C$ is again toroidal and weakly
semistable; in particular, $X_C$ has canonical singularities.
\end{enumerate}
\end{lem}

{\bf Proof.} Part 1. is proved in \cite{ak}, Lemma~6.1. Part 2. is proved
in \cite{ak}, Lemma~6.2 for the case when $g:C\rar B$ is a dominant
morphism; the proof for a curve $C$ is the same word-by-word. \QED

\begin{lem}\label{lem-wss-stab} Let $f: X\rar C$ be a weakly semistable
morphism from a variety $X$ with $\dimension X = n+1$ to a curve $C$. Assume
that the fibers of $f$ are of general type, and let $f_{can}:
X_{can}\rar C$ be the relative canonical model. Then the special fiber
$X_{can,0}=f^{-1}_{can}(0)$ is a stable n-fold. 
\end{lem}

{\bf Proof.} By Lemma~\ref{lem-wss-prop}, if $g: (C',0')\rar (C,0)$ is a finite
base change, then $X'=X\times_C C'$ is again canonical. Since $K_{X'} = g^*
K_X + m X'_{0'}$, we see that finding the canonical model commutes with taking a
finite base change:
\[ X_{can}\times_C C' \cong (X\times_C C')_{can} \]
and we can apply Theorem~\ref{lem-ksb}~(ii) to conclude that $X_{can,0}$  has
semi-log-canonical singularities. \QED 

It can be shown directly that the fibers of a
weakly semistable morphism have semi-log-canonical singularities. We may
assume that the base is a curve. Then in a local model the special fiber is
the complement of the big torus in a toric variety. Alexeev has shown
\cite{ale2} that this complement is semi-log-canonical.

\begin{lem}\label{lem-gg} Let $f:X\rar C$ be as in the previous lemma, and assume that the
genus $g(C)\geq 3$. Then the sheaf $f_* \co_X (l K_X)$ is generated by global
sections for all $l\geq 2$. 
\end{lem}

{\bf Proof.} (Koll\'ar \cite{kol}) By Theorem 4.12 in \cite{kol}, the vector
bundle $E_l=f_* \co_X(l K_{X/C})$ is semipositive for all $l\geq 0$; that
means, any quotient line bundle of $E_l$ has non-negative degree on $C$. If
$L$ is a line bundle with $\mbox{deg}(L) > 2 g(C)-2$ then 
$\co_C(K_C)$ cannot be a quotient of $E_l\otimes L$. Hence $H^1(C,E_l\otimes
L)=0$. Now if $\mbox{deg}(L) > 2g(C)$ then $E_l\otimes L$ is generated by
global sections. In particular, $f_* \co_X (l K_X)= E_l \otimes \co_C(l K_C)$
is generated by global sections for all $l\geq 2$. \QED

\subsection{Deformations of canonical singularities and invariance of
plurigenera} 
The following two theorems were proved by Kawamata \cite{kaw}, generalizing
results of Siu \cite{siu}.

\begin{th}\label{thm-def-can} Let $\pi:X\rar S$ be a flat morphism from a germ
of a variety 
$(X,x_0)$ to a germ of a smooth curve $(S,s_0)$ whose special fiber $X_0 =
\pi^{-1}(s_0)$ has only canonical singularities. Then $X$ has only canonical
singularities. In particular, the fibers $X_s=\pi^{-1}(s)$ have only canonical
singularities. \qed
\end{th}

\begin{th}\label{thm-inv} Let $\pi:X\rar S$ be a projective flat morphism from
a normal variety to a germ of a smooth curve $(S,s_0)$. Assume that the
fibers $X_s = \pi^{-1}(s)$ have only canonical singularities and are of
general type for all $s\in S$. Then the plurigenus $P_m(X_s) = \dimension
H^0(X_s,m K_{X_s})$ is constant as a function on $s\in S$ for any $m>0$. \qed
\end{th} 

The two theorems for families of n-folds follow easily from the $MMP(n+1)$
asssumption. However, we are going to apply them 
for families of $(n+1)$-dimensional varieties, and we do not want to assume
minimal model program in dimension $n+2$.

\section {Proof of boundedness}

By Matsusaka's theorem (Theorem 2.4, \cite{mat2}) there exists an integer
$\nu_0>0$ such that 
if $X$ is a normal variety with rational Gorenstein singularities, with
ample canonical divisor $K_X$, and with given Hilbert function $H(l)=\chi(X,l
K_X)$, then 
$\nu_0 K_X$ is very ample and has no higher cohomology. Thus, the Hilbert
scheme 
parameterizing embeddings $X\in\bfp^{H(\nu_0)-1}$ with Hilbert function
$H(\nu_0 t)$ has 
a locally closed subscheme parameterizing $\nu_0$-canonical embeddings of
varieties $X$ with rational Gorenstein singularities. We give this subscheme
the induced reduced structure and call its closure $B_0$. Let $f_0: X_0\rar B_0$ be
the (closure of the) universal family over $B_0$. We apply weak semistable
reduction (Theorem~\ref{thm-wssr}) to the morphism $f_0$ to get a weakly
semistable morphism  
\[ f: X\rar B. \]  

By Lemma~\ref{lem-wss-prop} we know that for a general nonsingular curve
$C\subset B$ the inverse image $X_C=f^{-1}(C)$ has canonical
singularities. We can then apply the MMP(n+1) assumption and find a canonical
model $X_{C,can}$ for $X_C$. Since the general fiber of $X_{C,can}\rar C$ has
rational Gorenstein singularities, we see that the 
special fiber  is a stable smoothable n-fold with Hilbert
function $H$.

The idea of the proof is to construct a relative canonical model for the
whole family $f:X\rar B$ and to show that the restriction of this canonical
model to a curve $C$ is the canonical model $X_{C,can}$.

\begin{lem}\label{lem-fg} The relative canonical ring 
\[ R_{X/B} = \bigoplus_{l\geq 0} f_* \co_X (l K_X) \]
is a finitely generated $\co_B$-algebra.
\end{lem}

To prove the lemma we sweep $B$ locally with nonsingular curves. By Bertini's
theorem a general hyperplane section through $b\in B$ is nonsingular, and the
same is true for nearby hyperplane sections. Applying this $\dimension B-1$
times we get the following diagram:
\[
\begin{CD}
X_1 = B_1\times_B X @>>> X\\
@VVV @VV{f}V\\
B_1 @>{\phi}>> B \\
@V{g}VV \\
S
\end{CD}
\]
where
\begin{enumerate}
\item $(g:B_1\rar S,\phi)$ is a family of nonsingular curves in $B$
parameterized by $S$, no curve lying in $B\setmin U_B$;
\item there exist $b_1\in B_1$, $\phi(b_1)= b$, and open neighborhoods
$b_1\in U_1$, $b\in U$ such that $\phi: U_1\rar U$ is an isomorphism;
\item we may assume that the hyperplane sections have high degree, so the
fibers of $g$ will have high genus ($\geq 3$).  
\end{enumerate}

Since Lemma~\ref{lem-fg} is local in $B$ we may replace $B$ by $B_1$ and $X$
by $X_1$. 

\begin{lem}\label{lem-iso} Consider the Cartesian diagram
\[
\begin{CD}
Y @>>> X\\
@V{f}VV @VV{f}V\\
C @>>> B \\
@V{g}VV @VV{g}V \\
\{s\} @>>> S
\end{CD}
\]
where $s\in S$ is a closed point, $C=g^{-1}(s)$, and $Y=f^{-1}(C)$. Then the
natural morphism 
\[ f_* \co_X(l K_X) \otimes_{\co_B} \co_C \rar f_*\co_Y(l K_Y)\]
is an isomorphism for every $s\in S$ and every $l\geq 2$.
\end{lem}

{\bf Proof.} By ``cohomology and base change'' \cite{har} it suffices to prove the
surjectivity of this morphism. Lemma~\ref{lem-gg} shows that the sheaf
$f_*\co_Y(l K_Y )$ is generated by global sections $H^0(Y,\co_Y(l K_Y))$ for
every $l\geq 2$. The family $g\circ f: X\rar S$ has fibers of general type
(by sub-additivity of Kodaira dimension \cite{vie1}, since the fibers of both $f$ and $g$
are of general type) and with canonical singularities (by
Lemma~\ref{lem-wss-prop}). Theorem~\ref{thm-inv} now implies that  a 
global section in $H^0(Y,\co_Y(l K_Y))$ can be lifted to a local section of
$(g\circ f)_* \co_X(l K_X)$, which gives a section of $f_* \co_X(l K_X)$.
\qed

\begin{rem}\label{rem-ample-open}
For the next proof we note that ampleness is an open condition: given a
flat projective morphism $f:X\rar B$ and a $\bfq$-Cartier divisor $D$ on $X$
such that $D|_{X_b}$ is ample for some $b\in B$, then $D$ is $f$-ample over a
neighborhood of $b$ in $B$. This follows from Kleiman's numerical criterion
for ampleness \cite{kle}. If $Z_1(X/B)$ is the free group generated by reduced
irreducible curves in $X$ mapping to a point in $B$, and $N_1(X/B)$ is
$Z_1(X/B)\otimes \bfr$ modulo numerical equivalence (see \cite{kmm} for
notations), then we have natural linear maps $Z_1(X_b)\rar N_1(X/B)$ and
$Z_1(X_\eta)\rar N_1(X/B)$ where $\eta$ is the generic point of $B$. Now the
second map factors through the first: for a curve over the generic point of
$B$ we take its closure and restrict to the special fiber. Also, the maps
restricted to effective cones factor. If $D$ is ample on $X_b$ then $D$ is
positive on the effective cone $\overline{NE}(X_b)\setmin\{0\}$, hence it is
positive on $\overline{NE}(X_\eta)\setmin\{0\}$, and so $D$ is ample on
$X_\eta$.  
\end{rem}

{\bf Proof of Lemma~\ref{lem-fg}.} The statement is local in $B$ , so we may
replace $B$ by a smaller open neighborhood of $b\in B$ if necessary. We use the
same notation as in Lemma~\ref{lem-iso}.

By the minimal model program assumption MMP(n+1) applied to $f: Y\rar C$, the sheaf $R_{Y/C}$ is a
finitely generated $\co_C$-algebra, and by  Lemma~\ref{lem-iso},
$R_{X/B}\otimes \co_C \cong R_{Y/C}$, at least in high degrees. Since the
graded pieces 
$(R_{X/B})_l = f_*\co_X(l K_X)$ are finite $\co_B$-modules we can apply
Nakayama's lemma and conclude that $R_{X/B}\otimes
\co_{B,b}$ is a finitely generated $\co_{B,b}$-algebra. Consider
\[ \proj (R_{X/B}\otimes\co_{B,b}) \rar \spec \co_{B,b}. \]
Since this morphism is projective, defined by finitely many
equations, we get a projective scheme $\tilde{X}$ over an open neighborhood
$U\subset B$ of $b$, and after replacing $B$ by $U$:
\[ \tilde{f}:\tilde{X}\rar B.\] 

By construction, $\tilde{Y}=\tilde{f}^{-1}(C)$ is a relative canonical model
of $f: Y\rar C$; in particular, it has canonical singularities. By
Theorem~\ref{thm-def-can} canonical singularities deform to 
canonical singularities, hence $\tilde{X}$ has canonical singularities. 

The canonical divisor $K_{\tilde{X}}$ is $\tilde{f}$-ample when restricted to
$\tilde{Y}$. Using the numerical criterion for ampleness,
$K_{\tilde{X}}$ is $\tilde{f}$-ample over a neighborhood of
$b$.  Again, we may replace $B$ by a smaller open set and assume that
$K_{\tilde{X}}$ is $\tilde{f}$-ample. 

Now $\tilde{f}:\tilde{X}\rar B$ is a
relative canonical model, birational to $f: X\rar B$ over $B$, hence the
relative canonical rings are the same
\[ R_{X/B} \cong R_{\tilde{X}/B}.\]
Since the latter is finitely generated as an $\co_B$-algebra, so is the
former. \QED

{\bf Proof of Theorem~\ref{thm-bound}}
We show that $\tilde{f}: \tilde{X} = \proj R_{X/B} \rar B$ is the
required family. 

It is clear from the construction that every closed fiber $\tilde{f}^{-1}(b)$
is a stable smoothable n-fold with Hilbert function $H$. Indeed, we have
constructed through every closed point $b\in B$ a nonsingular curve $C$ such
that $\tilde{f}^{-1}(C)$ is the relative canonical model of a weakly
semistable family over $C$ with
Gorenstein generic fiber, so the special fiber is stable and smoothable by
Lemma~\ref{lem-wss-stab}. 

To prove the converse, namely that the fibers of $\tilde{f}$ contain all  stable
smoothable n-folds with Hilbert function $H$, we go back
to the construction of $\tilde{f}:\tilde{X}\rar B$ from the universal family
over the Hilbert scheme. We have the commutative diagram 
\[
\begin{CD}
\tilde{X} @<<< X @>>> X_1 = X_0\times_B B_0 @>>> X_0\\
@V\tilde{f}VV  @V{f}VV  @V{f_1}VV  @V{f_0}VV\\
B @<{=}<< B @>{=}>> B @>>> B_0 \\
\end{CD}
\]
where the morphisms $B\rar B_0$ and $X\rar X_1$ come from weak semistable
reduction and the map $X\rar\tilde{X}$ is only rational.

Let $V_{B_0}\subset B_0$ be the open set parameterizing $\nu_0$-canonical
embeddings of normal 
varieties with rational Gorenstein singularities; and let $V_B\subset B, 
V_{X_1}\subset X_1$, and $V_{\tilde{X}}\subset \tilde{X}$ be the inverse
images of $V_{B_0}$. Now $f_1: 
V_{X_1}\rar V_B$ is a family of normal varieties with rational Gorenstein
singularities and ample canonical bundle over a nonsingular
base. It follows that $V_{X_1}$ is normal with canonical singularities and
$f_1: V_{X_1}\rar V_B$ is a relative canonical model, birational  to
$\tilde{f}:V_{\tilde{X}}\rar V_B$ over $V_B$. By uniqueness of the relative
canonical model, the two are isomorphic $V_{\tilde{X}}\cong V_{X_1}$ over
$V_B$. 

Let $\pi:Y\rar C$ be a smoothing of a stable n-fold
$Y_0$, and let $Y_\eta = \pi^{-1}(\eta)$ be the fiber over the generic point
$\eta$ of $C$. Since $Y_\eta$ is a projective rational Gorenstein scheme
with ample canonical bundle and Hilbert function $H$, we can find a morphism
$\eta\rar V_{B_0}$ such that the universal family $V_{X_0}$ pulls back to
$Y_\eta$. There exists a finite morphism $(C',\eta')\rar (C,\eta)$ such that
we can lift  $\eta\rar V_{B_0}$ to $\eta' \rar V_B$. By
completeness of $B$ we get a morphism $C'\rar B$. Let $\tilde{Y} = C'\times_B
\tilde{X}$, and let $Y'=C'\times_C Y$. The two families $\tilde{Y}\rar C'$
and $Y'\rar C'$ are isomorphic over the generic point, so their special
fibers must be isomorphic by separatedness of the functor $\cm_H^{sm}$
(Lemma~\ref{lem-sep}).  \qed

\subsection{Boundedness for stable smoothable pairs}
With the notations of \cite{ale2}, we consider stable pairs $(Y_0,D_0)$ where
$Y_0$ is a connected projective variety and $D_0$ is a reduced effective Weil
divisor in $Y_0$ such that the pair $(Y_0,D_0)$ has semi-log-canonical
singularities and $K_{Y_0}+D_0$ is an ample $\bfq$-Cartier divisor. We say
that $(Y_0,D_0)$ is smoothable if it admits a deformation to a stable pair
$(Y_t,D_t)$  where $Y_t$ has rational singularities and $K_{Y_t}+D_t$ is
Cartier. Further, we fix the Euler characteristics $\chi(Y_0,
l(K_{Y_0}+D_0))$ and $\chi(D_0,l(K_{Y_0}+D_0)|_{D_0})$ for all $l$.

By Matsusaka's theorem, for some $\nu_0>0$ the divisor $\nu_0(K_{Y_t}+D_t)$
is very ample and without higher cohomology for all pairs $(Y_t,D_t)$ as
above. There exists a Hilbert scheme parameterizing embeddings $Y_t\subset
\bfp^N$ with the given Hilbert function. Since $D_t \subset Y_t\subset
\bfp^N$ also has a fixed Hilbert function, there is a Hilbert scheme
parameterizing such embeddings. In the product of the two Hilbert schemes we
can find a locally closed subscheme parameterizing the required embeddings
$(Y_t,B_t)\subset \bfp^N$ via $\nu_0(K_{Y_t}+D_t)$; there is also a universal
family $(\cal{Y},\cal{D})$ over this subscheme. 

As before, we complete the universal family and apply weak semistable
reduction to it. We can make sure that the proper transform $D$ of the universal
divisor $\cal{D}$ in the weakly semistable family $f:X\rar B$ is the union of 
horizontal toroidal divisors. The proof that  $\bigoplus_l \co_X (l(K_X+D))$ is a
finitely generated $\co_B$-algebra goes as before, replacing the $MMP(n+1)$
with $log MMP(n+1)$ assumption, if Koll\'{a}r's
semipositivity result and Siu-Kawamata theorems on invariance of plurigenera
and deformations of canonical singularities also hold for
pairs. I do not know if these are true. Koll\'{a}r's
semipositivity theorem can be avoided in the proof of boundedness by
considering ramified cyclic covers of the base variety $S$; but to construct the
moduli space we would still need it. As remarked above, the theorems on
invariance of plurigenera and deformations of canonical singularities can be
replaced by the $MMP(n+2)$ assumption; this is, however, a very strong
assumption, considering for example the case of surfaces $n=2$.

\end{document}